\documentclass[twocolumn,10pt]{IEEEtran}
\usepackage{graphicx}
\usepackage{psfrag}
\usepackage[tight,footnotesize]{subfigure}
\usepackage{amsmath}
\usepackage{amssymb}
\usepackage{mathrsfs}
\usepackage{stfloats}
\usepackage{dsfont}
\usepackage{setspace}
\usepackage{array}
\usepackage{bbm}
\usepackage{mathabx}
\usepackage{algorithmic}
\usepackage{algorithm}
\usepackage{glossaries}
\usepackage{color}
\usepackage{bm}

\newcommand{\eeq}{\end{equation}}

\newcommand{\bq}{\mbox{\boldmath $q$}}

\newcommand{\bh}{\mbox{\boldmath $h$}}
\newcommand{\bH}{\mbox{\boldmath $H$}}

\newcommand{\bG}{\mbox{\boldmath $G$}}
\newcommand{\bPhi}{\mbox{\boldmath $\Phi$}}

\newcommand{\bu}{\mbox{\boldmath $u$}}

\newcommand{\bx}{\mbox{\boldmath $x$}}

\newcommand{\bg}{\mbox{\boldmath $g$}}

\newcommand{\bv}{\mbox{\boldmath $v$}}

\newcommand{\by}{\mbox{\boldmath $y$}}
\newcommand{\ds}{\displaystyle}
\newcommand{\bw}{\mbox{\boldmath $w$}}

\newcommand{\beq}{\begin{equation}}

\newtheorem{proposition}{Proposition}

\makeglossaries

\newacronym{mac}{MAC}{multiple-access channel}
\newacronym{bc}{BC}{broadcast channel}
\newacronym{mimo}{MIMO}{multiple-input multiple-output}
\newacronym{siso}{SISO}{single-input single-output}
\newacronym{sc}{SC}{single-carrier}
\newacronym{mc}{MC}{multi-carrier}
\newacronym{ofdma}{OFDMA}{orthogonal frequency division multiple access}
\newacronym{af}{AF}{amplify-and-forward}
\newacronym{df}{DF}{decode-and-forward}
\newacronym{cf}{CF}{compress-and-forward}
\newacronym{mwrc}{MWRC}{multi-way relay channel}
\newacronym{pde}{PDE}{partial data exchange}
\newacronym{fde}{FDE}{full data exchange}
\newacronym{iid}{i.i.d.\@}{independent and identically distributed}
\newacronym{awgn}{AWGN}{additive white Gaussian noise}
\newacronym{awg}{AWG}{additive white Gaussian}
\newacronym{sic}{SIC}{successive interference cancellation}
\newacronym{snr}{SNR}{signal-to-noise ratio}
\newacronym{sinr}{SINR}{signal-to-interference-plus-noise ratio}
\newacronym{ber}{BER}{bit error rate}
\newacronym{zf}{ZF}{zero-forcing}
\newacronym{mmse}{MMSE}{minimum mean square error}
\newacronym{sud}{SUD}{single user decoding}
\newacronym{dof}{DoF}{degrees of freedom}
\newacronym{gdof}{GDoF}{generalized degrees of freedom}
\newacronym{nnc}{NNC}{noisy network coding}
\newacronym{dmn}{DMN}{discrete memoryless network}
\newacronym{csi}{CSI}{channel state information}
\newacronym{ee}{EE}{energy efficiency}
\newacronym{ian}{IAN}{treating interference as noise}
\newacronym{snd}{SND}{simultaneous non-unique decoding}
\newacronym{brd}{BRD}{best response dynamics}
\newacronym{br}{BR}{best response}
\newacronym{ne}{NE}{Nash equilibrium}
\newacronym{lhs}{LHS}{left-hand side}
\newacronym{rhs}{RHS}{right-hand side}
\newacronym{gee}{GEE}{global energy efficiency}
\newacronym{wsee}{WSEE}{weighted sum energy efficiency}
\newacronym{wpee}{WPEE}{weighted product energy efficiency}
\newacronym{wmee}{WMEE}{weighted minimum energy efficiency}
\newacronym{kkt}{KKT}{Karush-Kuhn-Tucker}
\newacronym{pc}{PC}{pseudo-concave}
\newacronym{qc}{QC}{quasi-concave}
\newacronym{ql}{QL}{quasi-linear}
\newacronym{pl}{PL}{pseudo-linear}
\newacronym{spc}{SPC}{strictly pseudo-concave}
\newacronym{sqc}{SQC}{strictly quasi-concave}
\newacronym{lfp}{LFP}{linear fractional problem}
\newacronym{clfp}{CLFP}{concave-linear fractional problem}
\newacronym{ccfp}{CCFP}{concave-convex fractional problem}
\newacronym{mmfp}{MMFP}{max-min fractional problem}
\newacronym{sorp}{SoRP}{sum-of-ratios problem}
\newacronym{porp}{PoRP}{product-of-ratios problem}
\newacronym{srp}{SRP}{single-ratio problem}
\newacronym{brb}{BRB}{branch-reduce-and-bound}
\newacronym{qos}{QoS}{quality-of-service}
\newacronym{comp}{CoMP}{cooperative multi-point}
\newacronym{ue}{UE}{user equipment}
\newacronym{bs}{BS}{base station}
\newacronym{mrc}{MRC}{maximum ratio combining}
\newacronym{d2d}{D2D}{device-to-device}
\newacronym{lmmse}{LMMSE}{linear minimum mean square error}
\newacronym{ris}{RIS}{reconfigurable intelligent surface}
\newacronym{svd}{SVD}{singular values decomposition}

\begin{document}
\title{\huge Energy Efficiency Optimization of Reconfigurable Intelligent Surfaces with Electromagnetic Field Exposure Constraints}
\author{
	A. Zappone, {\em Senior Member, IEEE} and M. Di Renzo, {\em Fellow, IEEE}
\thanks{A. Zappone is with the University of Cassino and Southern Lazio, 03043 Cassino, Italy  (alessio.zappone@unicas.it). M. Di Renzo is with Universit\'e Paris-Saclay, CNRS, CentraleSup\'elec, Laboratoire des Signaux et Syst\`emes, 91192 Gif-sur-Yvette, France (marco.di-renzo@universite-paris-saclay.fr).}

\vspace{-0.75cm} }

\maketitle

\begin{abstract}
This work considers the problem of energy efficiency maximization in a \gls{ris}-based communication link, 
subject to not only the conventional maximum power constraints, but also additional constraints on the maximum exposure to electromagnetic radiations of the end-users. The \gls{ris} phase shifts, the transmit beamforming, the linear receive filter, and the transmit power are jointly optimized, and two provably convergent and low-complexity algorithms are developed. One algorithm can be applied to the general system setups, but does not guarantee global optimality. The second algorithm is provably optimal in a notable special case. The numerical results show that RIS-based communications can ensure high energy efficiency while fulfilling users' exposure constraints to radio frequency emissions.
\end{abstract}

\section{Introduction}
With fifth generation (5G) networks being rolled out, the attention is shifting towards the next generation of wireless networks. Among the candidate technologies for beyond 5G, the use of reconfigurable intelligent surfaces (RISs) is emerging as a way of ensuring high data rates, high energy efficiency, and a high degree of flexibility to adapt to sudden and heterogeneous service requests \cite{RIS_GE_JSAC,RuiZhang_COMMAG,SmartWireless}. Environmental objects separating the communication endpoints can be coated with \glspl{ris}, which can be dynamically reconfigured so as to provide channel customization features. In the area of radio resource allocation for \gls{ris}-based networks, most contributions focus on maximizing the network rate \cite{ZapTWC2021,Pan2019b}, the energy efficiency \cite{ZapTWC2021,EE_RISs}, or minimizing the power consumption \cite{Wu2018}.

At the same time, a relevant issue for future wireless networks is the growing concerns for electromagnetic pollution. Although, at present, non-ionizing radio frequency radiations have not been associated to any health condition \cite{ICNIRP}, the continuous exposure to electromagnetic fields (EMF) is a factor that raises concerns among end-users and diminishes their acceptance of emerging transmission technologies and massive network deployments \cite{Chiaraviglio2020}. A few studies have appeared in the literature to investigate resource allocation schemes that minimize the EMF exposure of human users. In \cite{Sambo15}, a survey on reducing the electromagnetic radiation in wireless systems is provided. In \cite{Wang2011}, a method for evaluating the specific absorption rate (SAR) of multi-antenna systems is proposed and evaluated. In \cite{Sambo2017}, the EMF exposure in orthogonal frequency-division multiplexing (OFDM) systems is minimized subject to minimum rate requirements for the users, and the authors of \cite{Castellano} have recently designed transmit policies that dynamically allocate users' electromagnetic radiation exposure over time.

While all these previous works considered the problem of EMF-aware communications in legacy wireless systems, an EMF-aware scheme for \gls{ris}-based wireless networks has recently been developed in \cite{Ibraiwish21}. Therein, the users' electromagnetic exposure is minimized subject to quality of service constraints. 
Instead, in this work we consider the different and more general problem of maximizing the energy efficiency in an \gls{ris}-assisted communication MIMO link, enforcing both maximum power constraints and maximum EMF exposure constraints. The optimization problem is tackled with respect to the \gls{ris} phase shifts, the transmit beamforming, the linear receive filter, and the transmit power. The EMF constraints are formulated in terms of maximum acceptable values for the SAR, which measures the rate of electromagnetic energy absorption per unit mass of human body when it is exposed to a radio frequency electromagnetic field \cite{Wang2011}, \cite{Castellano}.

In particular, we devise two provably convergent optimization methods for EMF-aware RIS-assisted communications. The first algorithm leverages alternating maximization and has a complexity that is linear in the number of \gls{ris} elements and polynomial in the number of transmit and receive antennas. The second algorithm achieves the global optimal solution without requiring any iterations, in the notable special case of isotropic EMF exposure constraints, enjoying a complexity that is linear in the number of \gls{ris} elements and in the number of transmit and receive antennas. The numerical results show that employing RISs can ensure high energy efficiency while fulfilling electromagnetic exposure constraints.

\section{System model}
Consider a single-user system in which a transmitter with $N_{T}$ antennas and a receiver with $N_{R}$ antennas communicate through an \gls{ris}. The direct link between the transmitter and receiver is assumed to be weak enough to be ignored. Denote by $\delta$ the end-to-end path loss, $\bH$ and $\bG$ the fading channels from the transmitter to the \gls{ris} and from the \gls{ris} to the receiver, respectively, $p$ the transmit power, $\bq$ and $\bw$ the unit-norm transmit beamformer and receive combiner. The \gls{ris} has $N$ elementary passive scatterers, which can independently reflect the radio wave impinging upon them according to a unit amplitude reflection coefficient $e^{j\phi_{n}}$, $n=1,\ldots,N$ and $j$ denoting the imaginary unit. We assume that a reliable channel estimation phase has been performed, and that the \gls{ris} phase configuration can be set by sending a configuration signal to an \gls{ris} controller with minimal signal processing, transmission/reception, and power storage capabilities \cite[Fig. 4]{RIS_GE_JSAC},  \cite{ZapTWC2021}. Under these assumptions, the system bit-per-Joule energy efficiency is expressed as 
\begin{align}\label{Eq:EE}
\text{EE}&=\frac{B\log_2\left(1+\ds\frac{p}{\ds\delta \sigma^{2}}\left|\bw^{H}\bG\bPhi\bH\bq\right|^{2}\right)}{\mu p+P_{c}}
\end{align}
with $\mu$ the inverse of the transmit amplifier efficiency, $P_{c}$ the static power consumption of the system, $B$ the communication bandwidth, and $\sigma^{2}$ the receive noise power. 

\subsection{EMF-Aware Optimization: Near-Field SAR Constraints}
Each portable device must comply with specific SAR limits of radiation that are considered safe for the body \cite{Castellano}. The challenge of fulfilling SAR compliance is exacerbated by the use of multiple transmit antennas in portable wireless devices, which increases the exposure for a given total transmit power, due to the combinations of the precoding gains and phases across the antennas \cite{Wang2011}. The objective of this paper is to optimize the \gls{ris} matrix $\bPhi$, the beamforming vector $\bq$,  the receive filter $\bw$, and the transmit power $p$, for energy efficiency maximization, subject to both power and EMF constraints. 

To this end, we formulate the following problem
\begin{subequations}\label{Prob:EMF}
\begin{align}
&\ds\max_{\bPhi,\bq,\bw}\;\frac{B\log_2\left(1+\ds\frac{p}{\ds\delta \sigma^{2}}\left|\bw^{H}\bG\bPhi\bH\bq\right|^{2}\right)}{\mu p+P_{c}}\label{Prob:aEMF}\\
&\;\text{s.t.}\;\phi_{n}\in[0,2\pi]\;,\;0\leq p\leq P_{max}\label{Prob:aEMF}\\
&\;\quad\;\sum\nolimits_{n=1}^{N_{T}}c_{n}|q_{n}|\leq P_{q}\;,\;\sum\nolimits_{n=1}^{N_{T}}|q_{n}|^{2}\leq 1\label{Prob:cEMF}\\
&\;\quad\;\sum\nolimits_{n=1}^{N_{R}}d_{n}|w_{n}|\leq P_{w}\;,\;\sum\nolimits_{n=1}^{N_{R}}|w_{n}|^{2}\leq 1\;,\label{Prob:dEMF}
\end{align}
\end{subequations}
wherein $P_{max}$ is the maximum transmit power, $P_{q}$ and $P_{w}$ are the maximum EMF constraints at the transmitter and receiver side, respectively, while $\{c_{n}\}_{n=1}^{N}$ and $\{d_{n}\}_{n=1}^{N}$ are the EMF absorption coefficients, which account for the magnitude of the total electric field absorbed by the human body due to the beamformer applied by the $n$-th transmit antenna, with $n=1,\ldots,N_{T}$, and the receive combiner applied by the $n$-th receive antenna, with $n=1,\ldots,N_{R}$. The EMF constraint does not apply to the RIS because of the unit modulus design assumption of the matrix of reflection coefficients.

\section{Energy efficiency maximization}\label{Sec:Optimization}
The EE maximization problem will be first tackled in its general form given in \eqref{Prob:EMF}, for which an iterative method will be developed. Next, a special case of Problem \eqref{Prob:EMF} will be globally solved in closed-form. 

\subsection{Maximization by Alternating Optimization}\label{Sec:AO}
A suitable approach to tackle Problem \eqref{Prob:EMF} is the alternating optimization of the \gls{ris} phase shift matrix $\bPhi$, the beamforming vector $\bq$, the receive filter $\bw$, and the transmit power $p$. These four subproblems are solved in the next three sections. 
\subsubsection{Optimal $\bPhi$}\label{Sec:OptPhi}
For fixed $\bq$, $\bw$, $p$, the problem becomes 
\begin{subequations}\label{Prob:EMF_Phi}
\begin{align}
&\ds\max_{\bPhi}\;|\bw^{H}\bG\bPhi\bH\bq|\label{Prob:aEMF_Phi}\;,\;\text{s.t.}\;\phi_{n}\in[0,2\pi],
\end{align}
\end{subequations}
since the denominator does not depend on $\bPhi$ and the logarithm is an increasing function. Problem \eqref{Prob:EMF_Phi} is solved by setting $\phi_{n}=-\angle{g_{n}^{*}h_{n}}$, where $(\cdot)^{*}$ is the complex conjugate, $g_n$ and $h_n$ are the $n$-th component of $\bG^{H}\bw$ and $\bH\bq$, respectively. 
\subsubsection{Optimal $\bq$}\label{Sec:Optq}
For fixed $\bPhi$, $\bw$, $p$, the optimization becomes 
\begin{align}\label{Prob:EMF_q}
&\ds\max_{\bq}\;|\bv^{H}\bq|\;,\;\text{s.t.}\sum\nolimits_{n=1}^{N_{T}}c_{n}|q_{n}|\leq P_{q}\;,\sum\nolimits_{n=1}^{N_{T}}|q_{n}|^{2}\leq 1\;,
\end{align}
wherein $\bv^{H}=\bw^{H}\bG\bPhi\bH$, and $P_{q}$ is the maximum allowed EMF exposure due to the transmit antennas. Since the constraints involve only the moduli of the components of $\bq$, it is optimal to set the phases of $q_{n}$ so as to align the phases of the entries $v_{n}^{*}$ in $\bv^{H}$. Plugging $\angle{q_{n}}=\angle{v_{n}}$ into \eqref{Prob:EMF_q} yields 
\begin{align}\label{Prob:EMF_q2}
&\ds\max_{\{x_{n}\geq 0\}_n}\sum_{n=1}^{N_{T}}|v_{n}|x_{n}\;,\text{s.t.}\sum_{n=1}^{N_{T}}c_{n}x_{n}\leq P_{q}\;,\sum_{n=1}^{N_{T}}x_{n}^{2}\leq 1
\end{align}
with $x_{n}=|q_{n}|$ for $n=1,\ldots,N$. Problem \eqref{Prob:EMF_q2} is convex and thus can be globally solved with polynomial complexity.
\subsubsection{Optimal $\bw$}\label{Sec:Optw}
Defining $\bu=\bG\bPhi\bH\bq$ and denoting by $P_{w}$ the maximum allowed EMF exposure due to the receive antennas, the same line of reasoning used to optimize $\bq$ leads us to setting $\angle{w_{n}}=\angle{u_{n}}$, which yields the problem  
\begin{align}\label{Prob:EMF_w}
&\ds\max_{\{y_{n}\geq0\}_n}\sum_{n=1}^{N_{R}}|u_{n}|y_{n}\;,\text{s.t.}\sum_{n=1}^{N_{R}}d_{n}y_{n}\leq P_{w}\;,\sum_{n=1}^{N_{R}}y_{n}^{2}\leq 1
\end{align}
with $y_{n}=|w_{n}|$ for $n=1,\ldots,N$. Problem \eqref{Prob:EMF_w} is convex and thus can be globally solved with polynomial complexity.
\subsubsection{Optimal $p$}\label{Sec:OptP} The optimal $p$ is found as 
\begin{align}\label{Prob:EMF_p}
&\ds\max_{p}\;\frac{\log_{2}(1+pc)}{\mu p+P_{c}}\;,\text{s.t.}\;0\leq p\leq P_{max}
\end{align}
with $c=|\bw^{H}\bG\bPhi\bH\bq|^{2}/\delta \sigma^{2}$. Problem \eqref{Prob:EMF_p} can be seen to be a pseudo-concave maximization, and thus its solution is simply found from the stationarity condition \cite{ZapNow15}.

\begin{algorithm}[!t]
\begin{algorithmic}\caption{Alternating maximization}
\label{Alg:AO}
\small
\STATE \texttt{Initialize $\bw$ and $\bq$ to feasible values.}
\REPEAT
\STATE $\bg=\bw^{H}\bG$; $\bh=\bH\bq$; $\phi_{n}=-\angle{g_{n}^{*}h_{n}}\;, \forall n=1,\ldots,N$;
\STATE $\bv^{H}=\bw^{H}\bG\bPhi\bH$; $q_{n}=x_{n}e^{j\angle{v_{n}}}\;,\forall n=1,\ldots,N$, \texttt{with } $\bx$ \texttt{ the solution of Problem \eqref{Prob:EMF_q2}};
\STATE $\bu=\bG\bPhi\bH\bq$; $w_{n}=y_{n}e^{j\angle{u_{n}}}\;, \forall n=1,\ldots,N$, \texttt{with } $\by$ \texttt{ the solution of Problem \eqref{Prob:EMF_w}}; 
\UNTIL{Convergence}
\STATE \texttt{Set }$p$ \texttt{ as the solution of \eqref{Prob:EMF_p}};
\end{algorithmic}
\end{algorithm} \setlength{\textfloatsep}{5pt}%

Finally, the overall optimization algorithm can be stated as in Algorithm \ref{Alg:AO}. Let us observe that, the optimization of $p$ needs to be performed only once, after convergence has been reached for $\bq$, $\bw$, $\bPhi$, since the optimization of $\bq$, $\bw$, $\bPhi$ does not depend on the optimal $p$. As for the convergence of Algorithm \ref{Alg:AO}, the following result holds. 
\begin{proposition}
Algorithm \ref{Alg:AO} monotonically improves the value of the energy efficiency and converges.
\end{proposition}
\begin{IEEEproof}
Each step of Algorithm \ref{Alg:AO} globally optimizes one of the problem variables, thus leading to an increase of the objective \eqref{Prob:aEMF}. Therefore, the energy efficiency monotonically increases after each step of Algorithm \ref{Alg:AO}. Since the energy efficiency is an upper-bounded function, Algorithm \ref{Alg:AO} must converge in the value of the objective.  
\end{IEEEproof}

\subsection{A special case of Problem \eqref{Prob:EMF}}
Aiming at providing a deeper insight into the structure of the solution of Problem \eqref{Prob:EMF}, let us consider the notable special case in which $c_{n}=c$ for $n=1,\ldots,N_{T}$ and $d_{n}=d$ for $n=1,\ldots,N_{R}$. This case study is of practical relevance because it is likely that, for typical co-located multi-antenna devices, the human body has the same SAR for every antenna of the array. 
Next, we also assume that $P_{q}/c \leq 1$ and $P_{w}/d\leq 1$ hold, which corresponds to enforcing particularly demanding EMF constraints. Under these assumptions, the EMF constraints become $\sum_{n=1}^{N_{T}}|q_{n}|\leq {P_{q}}/{c}$ and $\sum_{n=1}^{N_{T}}|w_{n}|\leq {P_{w}}/{d}$. Moreover, they imply the unit-norm constraints, because it holds that 
\begin{align}
\sum\nolimits_{n=1}^{N_{T}}|q_{n}|^{2}&\leq \sum\nolimits_{n=1}^{N_{T}}|q_{n}|\leq {P_{q}}/{c}\leq 1\;,\\
\sum\nolimits_{n=1}^{N_{R}}|w_{n}|^{2}&\leq \sum\nolimits_{n=1}^{N_{R}}|w_{n}|\leq {P_{w}}/{d}\leq 1\;,
\end{align}
At this point, let us revisit the results in Secs. \ref{Sec:Optq}, \ref{Sec:Optw}. 
\subsubsection{Optimization of $\bq$}\label{Sec:Opt_qEq}
Defining $\bv^{H}=\bw^{H}\bG\bPhi\bH$ and noting that the optimal phases of the entries of $\bq$ are $\angle{q_{n}}=\angle{v_{n}}$, the problem reduces to 
\begin{subequations}\label{Prob:EMF_qEq}
\begin{align}
&\ds\max_{\{x_{n}\geq 0\}_n}\;\sum\nolimits_{n=1}^{N_{T}}|v_{n}|x_{n}\;,\text{s.t.}\;\sum\nolimits_{n=1}^{N_{T}}x_{n}\leq {P_{q}}/{c}\;,
\end{align}
\end{subequations}
with $x_{n}=|q_{n}|$ for $n=1,\ldots,N_{T}$. The following result holds.
\begin{proposition}\label{Prop:Optq}
Denote by $n_{q}\in\{1,2,\ldots,N_{T}\}$ the index such that $|v_{n_{q}}|\geq |v_{n}|$, for $n=1,\ldots,N_{T}$. Then, the optimal solution of Problem \eqref{Prob:EMF_qEq} is given by 
\beq\label{Eq:Opt_q}
\left\{
\begin{array}{ll}
x_{n}=\frac{P_{q}}{c}\;,\; \text{for } n=n_{q}\\
x_{n}=0\;,\; \text{for } n\neq n_{q}
\end{array}\right.
\eeq
\end{proposition}
\begin{IEEEproof}
Since $|v_{n_{q}}|\geq |v_{n}|$ for $n=1,\ldots,N_{T}$, there exist non-negative values $\epsilon_{1},\ldots \epsilon_{N_{T}}$ such that $|v_{n}|=|v_{n_{q}}|-\epsilon_{n}$, for $n=1,\ldots,N_{T}$. Then, it holds that 
\begin{align}
& \sum\nolimits_{n=1}^{N_{T}}|v_{n}|x_{n}\!=\!|v_{n_{q}}|x_{n_{q}}\!+\!\sum\nolimits_{n\neq n_{q}}^{N_{T}}(
|v_{n_{q}}|-\epsilon_{n}) x_{n}\notag\\
 \!& \hspace{0.5cm}=\!|v_{n_{q}}|\sum\nolimits_{n=1}^{N_{T}}x_{n}\!-\!\sum\nolimits_{n\neq n_{q}}^{N_{T}}\epsilon_{n}x_{n}\leq |v_{n_{q}}|{P_{q}}/{c}\;,\label{Eq:Proof_q}
\end{align}
where the last inequality follows because $\sum_{n=1}^{N_{T}}x_{n}\leq P_{q}/c$ and $\sum_{n\neq\bar{n}}^{N_{T}}\epsilon_{n}x_{n}\geq 0$. Finally, the proof follows because the solution in \eqref{Eq:Opt_q} achieves the upper-bound in \eqref{Eq:Proof_q}.
\end{IEEEproof}

\subsubsection{Optimization of $\bw$}\label{Sec:Opt_wEq}
Defining $\bu^{H}=\bG\bPhi\bH\bq$ and noting that the optimal phases of the entries of $\bw$ need  to  fulfill  the  identity $\angle{w_{n}}=\angle{u_{n}}$, the problem to be solved reduces to 
\begin{subequations}\label{Prob:EMF_wEq}
\begin{align}
&\ds\max_{\{y_{n}\geq 0\}_n}\;\sum\nolimits_{n=1}^{N_{T}}|u_{n}|y_{n}\;,\text{s.t.}\sum\nolimits_{n=1}^{N_{T}}y_{n}\leq {P_{w}}/{d}\;,
\end{align}
\end{subequations}
with $y_{n}=|w_{n}|$ for $n=1,\ldots,N_{T}$. The following result holds.
\begin{proposition}\label{Prop:Optw}
Denote by $n_{w}\in\{1,2,\ldots,N_{T}\}$ the index such that $|u_{n_{w}}|\geq |u_{n}|$, for $n=1,\ldots,N_{T}$. Then, the optimal solution of Problem \eqref{Prob:EMF_wEq} is given by 
\beq\label{Eq:Opt_w}
\left\{
\begin{array}{ll}
y_{n}=\frac{P_{w}}{d}\;,\; \text{for } n=n_{w}\\
y_{n}=0\;,\; \text{for } n\neq n_{w}
\end{array}\right.
\eeq
\end{proposition}
The proof follows along the same line of reasoning as for the proof of Proposition \ref{Prop:Optq} and it is hence omitted.
Based on Propositions \ref{Prop:Optq} and \ref{Prop:Optw}, it is possible to obtain the optimal solution of Problem \eqref{Prob:EMF} when $\bq$, $\bw$, $\bPhi$, $p$ are jointly optimized. By direct inspection of Proposition \ref{Prop:Optq} and Proposition \ref{Prop:Optw}, in fact, we evince that the optimal $\bq$ and $\bw$ have only one non-zero component, and so they can be written as
\begin{align}\label{Eq:Opt_q1}
\bq&=\frac{P_{q}}{c}e^{j\phi_{q}\footnotesize(\bPhi,\bw)}[\underbrace{0,\ldots, 0}_{n_{q}-1},1,0,\ldots,0]\\
\label{Eq:Opt_w1}
\bw&=\frac{P_{w}}{d}e^{j\phi_{w}\footnotesize(\bPhi,\bq)}[\underbrace{0,\ldots, 0}_{n_{w}-1},1,0,\ldots,0]\;,
\end{align}
where we have highlighted the fact that the optimal $\phi_{q}$ will depend on $\bPhi$ and $\bw$ through the vector $\bv$, and that the optimal $\phi_{w}$ will depend on $\bPhi$ and $\bq$, through the vector $\bu$. In practice, \eqref{Eq:Opt_q1} and \eqref{Eq:Opt_w1} imply that, under the assumption of isotropic EMF constraints, the optimal solutions for the transmit beamforming and the receive decoding vectors consist of activating a single antenna at the transmitter and a single antenna at the receiver, which need to be appropriately chosen. Exploiting \eqref{Eq:Opt_q1} and \eqref{Eq:Opt_w1} we obtain
\begin{align}
&|\bw^{H}\bG\bPhi\bH\bq|=\left|\frac{P_{q}}{c}\frac{P_{w}}{d}e^{j\phi_{q}\footnotesize(\bPhi,\bw)}e^{-j\phi_{w}\footnotesize(\bPhi,\bq)}\bg_{n_{w}}^{T}\bPhi\bh_{n_{q}}\right|=\notag\\
\label{Eq:OptJoint}
&\frac{P_{q}}{c}\frac{P_{w}}{d}\left|\bg_{n_{w}}^{T}\bPhi\bh_{n_{q}}\right|\leq \frac{P_{q}}{c}\frac{P_{w}}{d}\sum_{n=1}^{N}|\bg_{n_{w}}(n)\bh_{n_{q}}(n)|\;,
\end{align}
wherein $\bg_{n_{w}}^{T}$ is the $n_{w}$-th row of $\bG$, $\bh_{n_{q}}$ is the $n_{q}$-th column of $\bH$, and the last inequality is obtained with equality upon optimizing $\bPhi$, i.e., by choosing $\phi_{n}=-\angle{g_{n_{w}}(n)h_{n_{q}}(n)}$ for $n=1,\ldots,N$, with $g_{n_{w}}(n)$ and $h_{n_{q}}(n)$ being the $n$-th components of $\bg_{n_{w}}$ and $\bh_{n_{q}}$. 

Thus, the phases $\phi_{q}(\bPhi,\bw)$ and $\phi_{w}(\bPhi,\bq)$ do not affect the value of the objective function, regardless of the values of $\bq$, $\bw$, $\bPhi$. Therefore, in order to determine the optimal $\bPhi$, $\bq$, and $\bw$, it remains to optimize the indexes $n_{w}$ and $n_{q}$. This is  difficult to be performed in closed-form, but it can be carried out through an exhaustive search over all possible $N_{T}N_{R}$ choices of the pair $(n_{w},n_{q})$. Finally, the optimal transmit power can be determined as in Section \ref{Sec:OptP} once the optimal $\bPhi$, $\bq$, and $\bw$ have been computed. 

Thus, the overall globally optimal resource allocation procedure can be given as in Algorithm \ref{Alg:Opt}. It computes \eqref{Eq:OptJoint} for each choice $(n_{w},n_{q})$, selects the choice $(n_{w}^{o},n_{q}^{o})$ that yields the largest value of \eqref{Eq:OptJoint}, and then allocates the system resources accordingly. 

\begin{algorithm}[!t]
\begin{algorithmic}\caption{Global optimization with $P_{q}/c\leq 1$, $P_{w}/d\leq 1$.}
\label{Alg:Opt}
\small
\FOR{$i=1$ \TO $N_{T}$}
\FOR{$j=1$ \TO $N_{R}$}
\STATE $n_{q}=i$; $n_{w}=j$; $\text{Obj}(i,j)=\sum_{n=1}^{N}|\bg_{n_{w}}(n)\bh_{n_{q}}(n)|$;
\ENDFOR
\ENDFOR
\STATE $(n_{q}^{o},n_{w}^{o})=\text{argmax}\;\text{Obj}(i,j)$;
\STATE $\bq\!=\!\frac{P_{q}}{c}[\underbrace{0,\ldots, 0}_{n_{q}^{o}-1},1,0,\ldots,0]$; $\bw\!=\!\frac{P_{w}}{d}[\underbrace{0,\ldots, 0}_{n_{w}^{o}-1},1,0,\ldots,0]$;
\STATE $\phi_{n}\!=\!-\angle{g_{n_{w}^{o}}(n)h_{n_{q}^{o}}(n)}$; \texttt{Set }$p$\texttt{ as the solution of \eqref{Prob:EMF_p}};
\end{algorithmic}
\end{algorithm}
\subsection{Complexity analysis}
Let us compare the complexity of Algorithms \ref{Alg:AO} and \ref{Alg:Opt}. 

As for Algorithm \ref{Alg:AO}, each of the subproblems can be easily solved, since $\bPhi$ can be optimized with linear complexity by computing $N$ times the quantity $\phi_{n}=-\angle{g_{n}^{*}h_{n}}$,  while $\bq$ and $\bw$ have a polynomial complexity in $N_{T}$ and $N_{R}$, respectively, since  Problems \eqref{Prob:EMF_q2} and \eqref{Prob:EMF_w} are concave\footnote{The order $\alpha$ of the polynomial is not available in closed-form, but a known worst-case bound is $\alpha=4$}. As for the complexity related to computing $p$, it is negligible since \eqref{Prob:EMF_p} is a scalar problem which needs to be solved just once. Finally, the overall complexity of Algorithm \ref{Alg:AO} scales linearly with the number $I$ of iterations required to reach convergence. 

Algorithm \ref{Alg:Opt} has a lower complexity, since it is not iterative (i.e. $I=1$), and the complexity is linear in $N$, $N_{R}$, and $N_{T}$. Indeed, the variable Obj$(i,j)$ needs to be computed $N_{T}N_{R}$ times, and each computation requires $N$ complex multiplications. Finally, the same negligible complexity as for Algorithm \ref{Alg:AO} is required to optimize $p$.

\section{Numerical Results}\label{Sec:Numerics}
In our numerical analysis, we set $B=5\,\textrm{MHz}$, $\delta=110\,\textrm{dB}$, $N_{0}=-174\; \textrm{dBm/Hz}$, $N_{T}=N_{R}=4$, $P_{c}=30\,\textrm{W}$, $P_{max}=20\,\textrm{W}$.  
As for the fading channels, a Rician model is considered, wherein $h_{n}\sim{\cal CN}(v_{h},1)$ and $g_{n}\sim{\cal CN}(v_{g},1)$, with $v_{h}$ and $v_{g}$ such that the power of the line-of-sight path is four times larger than the power of all the other paths. All results are averaged over $10^{3}$ independent channel realizations. Moreover, $c_{n}=c=1/N_{T}$ and $d_{n}=d=1/N_{R}$. For simplicity, thus, we consider the isotropic EMF setup.

\begin{figure}[!t]
\centering
\includegraphics[width=0.8\columnwidth]{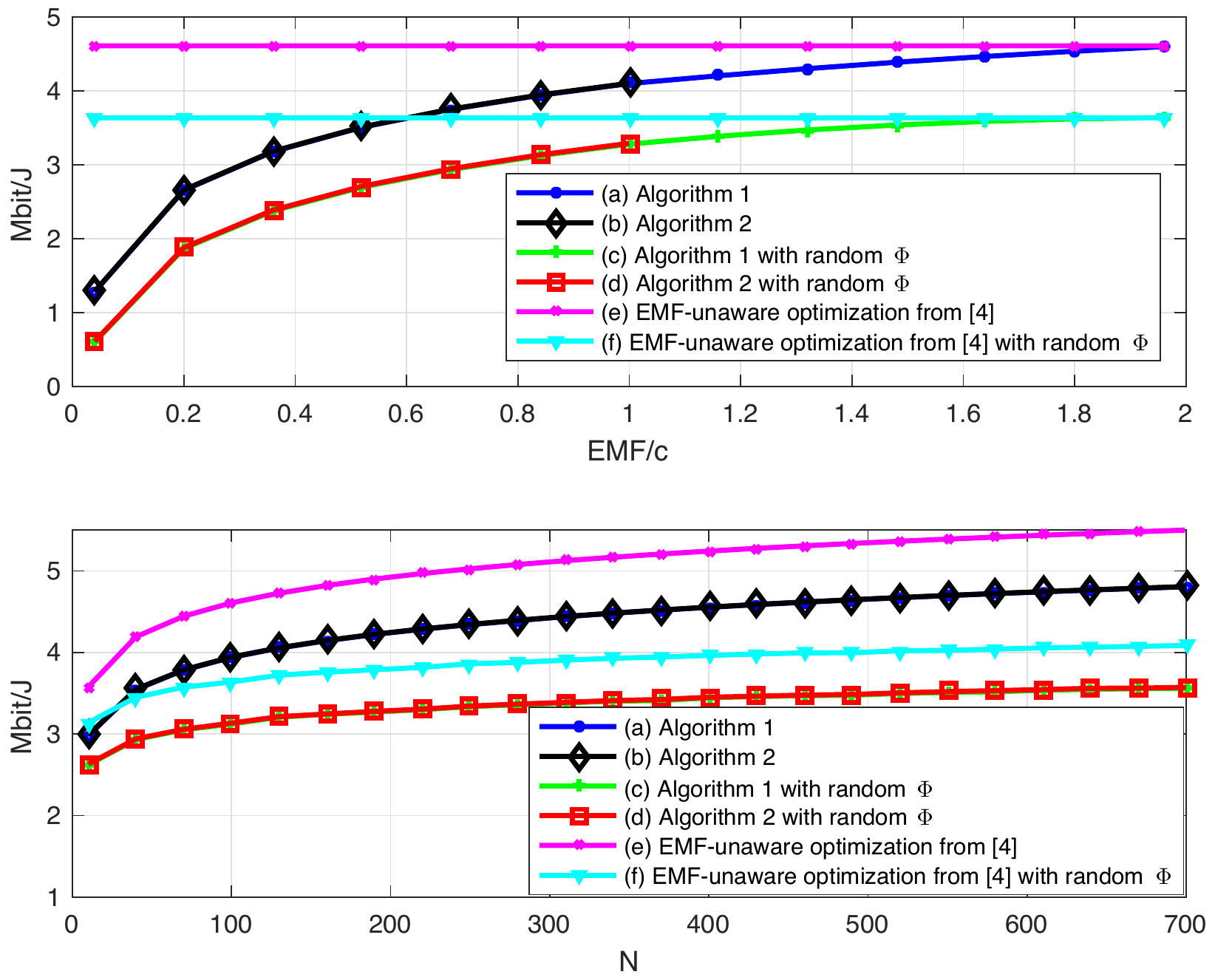}
\vspace{-0.5cm} \caption{Average energy efficiency of the six considered schemes as a function of (top) $P_{q}/c$ for $N=100$ and (bottom) $N$ for $P_{q}/c=0.85$.}
\label{Fig:EE}
\end{figure}

\begin{figure}[!t]
\centering
\includegraphics[width=0.9\columnwidth]{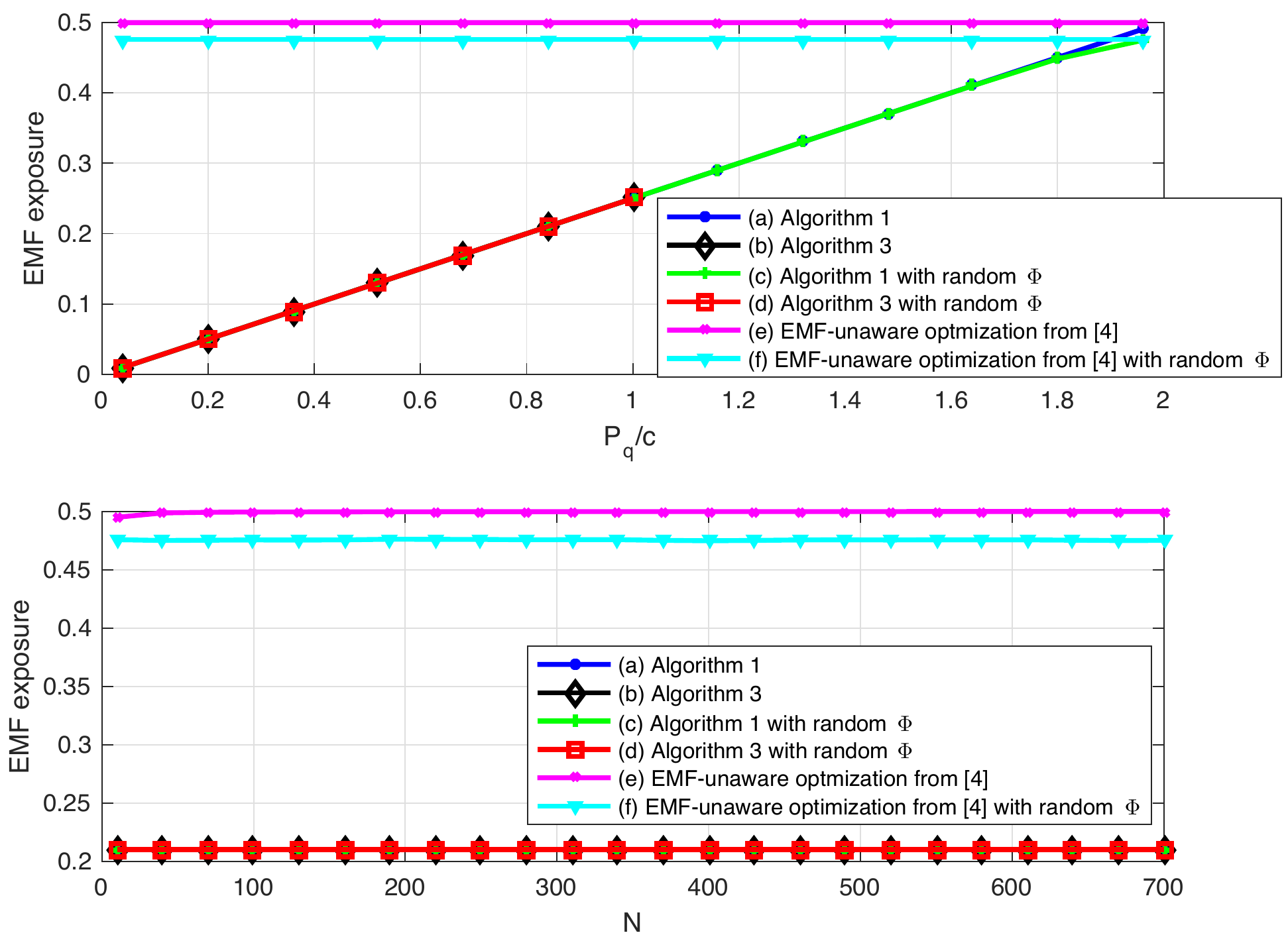}
\vspace{-0.5cm} \caption{EMF constraint $c\sum_{n=1}^{N}|q_{n}|$ of the six considered schemes as a function of (top) $P_{q}/c$ for $N=100$ and (bottom) $N$ for $P_{q}/c=0.85$.}
\label{fig:EMF}
\end{figure}

Figures \ref{Fig:EE} and \ref{fig:EMF} show the system energy efficiency and the EMF exposure\footnote{We focus on the exposure due to the transmit antennas. Similar results hold for the EMF exposure due to the receive antennas.} $c\sum_{n=1}^{N}|q_{n}|$, respectively, as a function of $P_{q}/c$ for $N=100$, and as a function of $N$ for $P_{q}/c=0.85$. The following schemes are considered and evaluated:
\begin{itemize}
\item[(a)] The EMF-aware alternating optimization of $\bq$, $\bw$, $\bPhi$, and $p$ by Algorithm \ref{Alg:AO}.
\item[(b)] The EMF-aware globally optimal optimization of $\bq$, $\bw$, $\bPhi$, and $p$ by Algorithm \ref{Alg:Opt}. Here, the curve is shown only in the range $P_{q}/c\leq 1$, since beyond this value the assumptions in Section \ref{Sec:Opt_qEq} do not hold. 
\item[(c)] The EMF-aware alternating optimization of $\bq$ and $\bw$ by Algorithm \ref{Alg:AO}. Here, the \gls{ris} matrix $\bPhi$ is not optimized and each phase shift is randomly set in $[0,2\pi]$. 
\item[(d)] The EMF-aware globally optimal optimization of $\bq$, $\bw$, and $p$, by Algorithm \ref{Alg:Opt}. Also in this case, each \gls{ris} phase shift is randomly set in $[0,2\pi]$. Moreover, this curve is shown in the range $P_{q}/c\leq 1$, since beyond this value the assumptions in Section \ref{Sec:Opt_qEq} do not hold. 
\item[(e)] The EMF-unaware optimization of $\bq$, $\bw$, $\bPhi$, and $p$ by the alternating optimization method from \cite{ZapTWC2021}. 
\item[(f)] The EMF-unaware optimization of $\bq$, $\bw$, and $p$, by the alternating optimization method from \cite{ZapTWC2021}. In this case, each \gls{ris} phase shift is randomly set in $[0,2\pi]$.
\end{itemize}

From Figs. \ref{Fig:EE} and \ref{fig:EMF}, we observe that enforcing an EMF constraint on the SAR of the human body reduces the energy efficiency level, since it restricts the feasible set of the problem. The use of RISs offers, however, the opportunity of achieving the desired energy efficiency while ensuring SAR-compliant communications. Figure \ref{Fig:EE} shows, in particular, that, by increasing $N$, we can attain the same energy efficiency as the benchmark systems in the absence of EMF constraints. EMF-aware transmission schemes lead, on the other hand, to large values of the EMF exposure $c\sum_{n=1}^{N}|q_{n}|$ (i.e., the EMF exposure), which may not fulfill the desired SAR values specified by national and international regulations. It is particularly interesting to note that increasing $N$ has little or no impact on the EMF exposure $c\sum_{n=1}^{N}|q_{n}|$ imposed to the transmit and receive filters. Also,  Algorithm \ref{Alg:AO} offers similar performance as the globally optimal solution obtained with Algorithm \ref{Alg:Opt}, although at a higher complexity. Moreover, as expected, the SAR-constrained energy efficiency tends, for large values of $P_{q}/c$, to the benchmark energy efficiency in the absence of EMF constraints. This is because increasing $P_{q}/c$ makes the EMF constraint become less relevant.

\section{Conclusions}
Low-complexity optimization algorithms have been proposed for energy efficiency maximization subject to EMF constraints. The analysis has shown that the use of a RIS can keep under control the end-users' EMF exposure while ensuring the desired energy efficiency level. Notably, this is obtained by using nearly-passive RISs that do not increase the amount of electromagnetic radiation over the air.

\bibliographystyle{IEEEtran}
\bibliography{FracProg}
\end{document}